\documentclass{article}

\begin{document}    
    
\input{amssym.def}    
\input{amssym.tex}    

 \newcommand{\R}{{\Bbb R}}    
 \newcommand{\Z}{{\Bbb Z}}    
 \newcommand{\C}{{\Bbb C}}    
 \newcommand{\N}{{\Bbb N}}    
\newcommand{\B}{{\cal B}}    
    
\title{Ordered abelian groups over a CW complex}

\author{Igor  ~Nikolaev\\
The Fields Institute\\
 222 College street\\
Toronto  M5T 3J1 Canada\\    
{\sf inikolae@fields.utoronto.ca}}

 \maketitle    
    
\newtheorem{thm}{Theorem}    
\newtheorem{lem}{Lemma}    
\newtheorem{dfn}{Definition}    
\newtheorem{rmk}{Remark}    
\newtheorem{cor}{Corollary}    
\newtheorem{prp}{Proposition}    
\newtheorem{con}{Conjecture}   
    
 \begin{abstract}    
If $X$ is a CW complex, one can assign to each point of $X$ an
ordered abelian group of finite rank whose subset of positive
elements depends continuously on the points of $X$. A locally
trivial bundle which arises in this way we denote by $E_X$. 
In the present work we establish a topological classification of such 
bundles in terms of the first cohomology group of $X$ with coefficients
in the ring ${\Bbb Z}_2$.

\vspace{7mm}    
    

\vspace{5mm}    
{\it AMS (MOS) Subj. Class.:  19K35, 46L40, 58F10}    
\end{abstract}

\section*{Introduction}    
In this note we prove the following theorem.

\bigskip\noindent
{\bf Theorem}\hskip0.3cm
{\it 
Let $X$ be a finite CW complex.
If $E_X$ is a bundle of the totally ordered abelian groups
of rank $k$ over $X$, then $E_X$ is classified by elements of
the first cohomology group $H^1(X;{\Bbb Z}_2)$. 
(See Section 2 for an exact formulation of this theorem.)
}

\bigskip\noindent
An interest to $E_X$ arises from the study of characteristic
classes of regular foliations on the hyperbolic manifolds,
see Section 4 of this note. Let us briefly remind this link.

To a noncommutative ring $R$ one can relate a semi-group
of equivalence classes of the finitely generated projective
modules over $R$. The completion of this semi-group to an abelian 
group is known as the Grothendieck, or $K_0$-group, of $R$.

To `restore' the ring $R$ from a $K_0$-group one is obliged
to introduce an extra structure on this group, called the {\it order}. 
Roughly speaking, the order on $K_0$ takes care of the initial semi-group  
of $K_0$. This order can be total or partial.

Let $X$ be a finite CW complex endowed with a foliation, $\cal F$.
The $C^*$-algebra of A.~Connes of $\cal F$ has  natural structure
of a noncommutative ring, $R$. The $K_0$-group of $R$ is isomorphic
to the first Betti group of $X$. This fact was established in
\cite{Nik1}, see also Ch. 10 of \cite{N}. The order on $K_0$
comes from the asymptotic behaviour of $\cal F$ at the universal
cover of $X$.

The homotopy class of $\cal F$ consists of foliations
obtained by continuous rotation of the tangent planes that
preserves the Frobenius (integrability) condition. There exist a parametrization
of the family of such foliations by the points of $X$ (Section 4).  
Thus, we have a bundle, $E_X$, of ordered abelian groups whose
`positive cone' depends continuously on the points $x\in X$.
The problem of classification of such bundles has an independent
interest and intrinsic beauty.

Theorem \ref{thm1} of Section 2 establishes a classification of $E_X$
in terms of the cohomology invariants of $X$. Surprisingly, these invariants
coincide with the elements of first cohomology group of $X$ with coefficients
in the ring ${\Bbb Z}_2$. 
One of applications of this theorem is calculation of the
characteristic classes of regular foliations (Corollary \ref{cor2}).
The characteristic classes have been extensively studied in the past 
by R.~Bott, A.~Haefliger, J.~Mather, W.~Thurston and others,
see a remarkable survey of Lawson \cite{Law}. This problem 
is known to be unsolved in the majority of cases, cf Haefliger  
p. 193 of \cite{Hae}.

\bigskip\noindent
The paper is organized as follows. In Section 1 we fix
notation and  terminology to be used throughout the paper.
We formulate and prove the main result in Sections 2 and 3, respectively.
An application of Theorem \ref{thm1} to the characteristic
classes of regular foliations is discussed in Section 4. 
An Appendix on the Stiefel-Whitney classes of vector bundles
is attached (Section 5). Finally, in Section 6 we suggest a generalization
of Theorem \ref{thm1} and Corollary \ref{cor2}.

\section{Notation}    
The partially ordered abelian groups make a category which is in a
sense dual to the category of noncommutative rings. This duality
appears if instead of usual isomorphisms between the rings one considers
rings which are Morita equivalent. The left (right) projective $R$-modules 
generate an abelian (Grothendieck) group with an order structure.
This order structure  defines a ring $R$ up to the Morita equivalence. 
An excellent introduction to the area is the monograph of Goodearl \cite{G}.

\subsection{Ordered abelian groups with interpolation}
Let $G$ be an additive abelian group. A {\it partial order} on $G$
is any reflexive, antysymmetric, transitive relation $\le$ on $G$.
If any pair $x,y\in G$ is comparable by this relation, the order
is called {\it total}. $G$ is a {\it partially ordered
abelian group} if for any $x,y,z\in G$ such that $x\le y$ it
follows that $x+z\le y+z$. The set of all positive elements 
$G^+\subset G$ is a {\it cone}, i.e. a subset of $G$ containing $0$ and 
closed under the addition. If $G$ is a totally ordered abelian
group, then
\begin{equation}\label{eq1}
G=G^+\bigcup (-G^+).
\end{equation}
An {\it order unit} in a partially ordered abelian group $G$
is a positive element $u\in G^+$ such that for any $x\in G$
there exists a positive integer $n$ such that $x\le nu$. Usually 
there are more than one choices of the order unit in a fixed group $G$.
The partially ordered abelian group $G$ with an order unit $u$
is denoted by $(G,u)$.    

\bigskip\noindent
\underline{Ideals of ordered groups}. 
Let $H\subseteq G$ be a subgroup of a partially ordered abelian group $G$.
$H$ is an {\it order ideal} in $G$ if for every $x,y\in H$ and $z\in G$
such that $x\le z\le y$, it follows that $z\in H$. The family of
order ideals in a partially ordered abelian group is a complete
lattice under the relation of inclusion (Corollary 1.10 of \cite{G}).
If $G$ has no order ideals except $\{0\}$ and $G$, then $G$ is
called {\it simple}.
\begin{lem}
For any order ideal $H\subseteq G$ the factor-group $G/H$    
is a partially ordered abelian group with the positive cone
\begin{equation}\label{eq2}
(G/H)^+\buildrel\rm def\over= (G^++H)/H.
\end{equation}
\end{lem}
{\it Proof.} See Goodearl \cite{G}.
$\square$

\bigskip\noindent
\underline{Interpolation and Elliott groups}.
The interpolation and Elliott groups represent two valuable classes of 
the partially ordered abelian groups. They were introduced
in connection with the study of linear operators in the Hilbert
space by F.~Riesz and classification of the AF $C^*$-algebras by 
G.~A.~Elliott, respectively.   

The partially ordered abelian group $G$ is called an {\it
interpolation group} if $G$ satisfies the Riesz interpolation
property: for any $x_1,x_2,y_1,y_2\in G$ such that $x_1\le y_1$,
$x_2\le y_1$, $x_1\le y_2$, $x_2\le y_2$, there exists $z\in G$
such that $x_1\le z\le y_1$, $x_2\le z\le y_1$, $x_1\le z\le y_2$
and $x_2\le z\le y_2$. In other words, there exist infinitely
many elements of the order structure lying `between' any two
elements. The interpolation property is hereditary with respect
to the order ideals $H$ and their quotients $G/H$ (Proposition 2.3
of \cite{G}).      

$G$ is called {\it unperforated} if for any positive integer
$n$ and the element $nx\ge 0$ it holds $x\ge 0$. By the
{\it dimension (Elliott) group} one understands an unperforated
partially ordered abelian group $G$ satisfying the Riesz interpolation
property. The order ideals and quotients of the Elliott groups are 
also the Elliott groups (Proposition 3.1 of \cite{G}). 

\bigskip\noindent
\underline{Simplicial groups.}
The simplicial groups are opposite to the simple partially ordered 
abelian groups. Such groups have an abundance of the order ideals.
On the other hand every Elliott group can be represented as the direct 
limit of a sequence of the simplicial groups.  

A {\it simplicial group} is a partially ordered abelian group ${\Bbb Z}^n$ 
whose positive cone consists of the vectors with nonnegative components:
\begin{equation}\label{eq3}
G^+=\sum_{i=1}^n {\Bbb Z}^+x_i,
\end{equation}
where $X=\{x_1,\dots,x_n\}$ is a basis in ${\Bbb Z}^n$.   
Every subset of $X$ generates an order ideal in the simplicial
group (Proposition 3.8 of \cite{G}). On the other hand, there
exist infinitely many possibilities to take a basis in ${\Bbb Z}^n$
so is the number of the order ideals in $G$. 
\begin{lem}
Any countable Elliott group is order isomorphic to a direct
limit of a countable sequence of simplicial groups.
\end{lem}
{\it Proof.} See Goodearl \cite{G}.
$\square$

\subsection{$S(G,u)$}
Let $(G,u)$ be a partially ordered abelian group with the order
unit. By a {\it state} on $(G,u)$ one understands a normalized
positive homomorphism from $(G,u)$ to the reals:
\begin{equation}\label{eq4}
s: G\longrightarrow {\Bbb R},\quad s(G^+)\subseteq {\Bbb R}^+,\quad
s(u)=1.
\end{equation}
The set of all states on $(G,u)$ is called a {\it state space} and
is denoted by $S(G,u)$. This is known to be a compact convex set
in the space ${\Bbb R}^G$ of all functions from $G$ to $\Bbb R$. 
As a convex set $S(G,u)$ is determined by the set of {\it extreme
points} and {\it faces} that correspond to certain states on $G$. 
Much of the theory of partially ordered abelian groups can be 
formulated in terms of the compact convex sets, cf Goodearl \cite{G}.

\bigskip\noindent
\underline{Discrete states}.
The homomorphism (\ref{eq4}) defines an additive subgroup in the
group of real numbers. This subgroup can be either a cyclic group
or a dense subset of $\Bbb R$ (Lemma 4.21 of \cite{G}). By a
{\it discrete state} one understands a state (\ref{eq4}) whose
$s(G)$ is a cyclic subgroup of $\Bbb R$. A link between the
discrete states and simplicial groups is given by the following lemma.
\begin{lem}
Let $s$ be a state on the interpolation group $(G,u)$. The set
$H=\{x-y~|~x,y\in (Ker ~s)^+\}$ is an order ideal of $G$. If
$s$ is a discrete state, then $G/H$ is a simplicial group. 
\end{lem}
{\it Proof.} See Goodearl \cite{G}.
$\square$ 

\bigskip\noindent
The extreme points of the space $S(G,u)$ can be  the
discrete states. By a {\it rational convex combination}
of the extreme points $x_1,\dots, x_n$ one understands a linear
combination $\alpha_1x_1+\dots+\alpha_nx_n$ such that $\alpha_i$
are non-negative rationals and $\sum_{i=1}^n\alpha_i=1$. 
Every discrete state on $(G,u)$ can be presented as a rational convex 
combination of the discrete extreme states (Proposition 6.22 of \cite{G}).

\bigskip\noindent
\underline{The space Aff $S(G,u)$}.
The state space $S(G,u)$ has been derived from a partially ordered
set. There exists a dual functor which `returns' $S(G,u)$ to a
partially ordered space, Aff $S(G,u)$, consisting of  affine
continuous real-valued functions on $S(G,u)$. This functor gives 
a useful `representation' of the group $(G,u)$.

The mapping $f: K_1\to K_2$ between two compact convex sets $K_1$ and 
$K_2$ is said to be {\it affine} if it preserves  the convex
combinations of points of $K_1$. Let $K_1=S(G,u)$ and $K_2={\Bbb R}$.
A {\it space Aff} $S(G,u)$ is defined as the space of all affine
continuous real-valued functions on $S(G,u)$. 

Let us define a map $\hat x: S(G,u)\to {\Bbb R}$ so that
$\hat x(s)=s(x)$ for all states $s\in S(G,u)$. There is
no difficulty to check that $\hat x$ is affine and continuous
mapping on $S(G,u)$. The function
\begin{equation}\label{eq5}
\phi: G\longrightarrow Aff ~S(G,u),
\end{equation}
acting by the rule $\phi(x)=\hat x$, is known as a
{\it natural map} from $G$ to Aff $S(G,u)$. (Note
that $\phi$ is not necessarily a surjective map.)
The following lemma gives a partial answer to a 
`representation problem' for the partially ordered abelian groups.
\begin{lem}\label{lm4}
Let $(G,u)$ be an interpolation group satisfying a general
comparability condition.
\footnote{This condition means that the state space of $(G,u)$
is `rich enough', see \cite{G} for an exact definition.
For example, every totally ordered abelian groups and simplicial
groups satisfy this condition.}
Let $\phi$ be the natural map of form (\ref{eq5}). Let us set
\displaymath
A=\{p\in Aff~S(G,u)~|~p(s)\in s(G)~\hbox{for all discrete} 
~s\in\partial_e S(G,u)\}, 
\enddisplaymath
where $\partial_e$ is the boundary of $S(G,u)$ consisting
of the extreme states. Then $\phi(G)$ is a dense subgroup
of $A$. 
\end{lem}
{\it Proof.} See Goodearl \cite{G}.
$\square$

\subsection{Ordered abelian groups of finite rank}
The results on ordered groups can be precised if one restricts to     
the groups whose state space is a Choquet simplex. The representation
theory for such groups (Lemma \ref{lm4}) is relatively full developed
area, see Goodearl \cite{G} for the details. In the case $G$ is
a simple Elliott group of finite rank, a complete classification is
possible.

\bigskip\noindent
\underline{Choquet simplex}.
A {\it Choquet simplex} is a compact simplex whose inductive dimension
is allowed to be infinity. For example, if $(G,u)$ is an interpolation
group, then $S(G,u)$ is a Choquet simplex (Theorem 10.17 of \cite{G}). 
The fundamental advantage is that any face of such a simplex is a
simplex (of lower dimension). In particular, every Choquet simplex
of finite dimension coincides with the usual simplex.

Until the end of this paragraph we assume that $(G,u)$ is simple 
(i.e. has only trivial order ideals). In this case any nonzero
element of $G^+$ can be chosen as the order unit $u$ (Lemma 14.1
of \cite{G}). The `majority' of simple groups cannot afford any
discrete states because of the following lemma.
\begin{lem}\label{lm5}
Let $(G,u)$ be a simple interpolation group. Then the following three 
conditions are equivalent: (i) there exists a discrete state on
$(G,u)$; (ii) $G$ is a cyclic abelian group; (iii) $(G,u)$ is
order isomorphic to $({\Bbb Z},m)$ for some $m\in {\Bbb N}$.
\end{lem}
{\it Proof.} See Goodearl \cite{G}.
$\square$

\bigskip\noindent
As a corollary of Lemma \ref{lm5} one obtains the representation theory of 
simple interpolation groups without perforation (Elliott groups). Namely,
the image $\phi(G)$ of the natural map (\ref{eq5}) is a dense subgroup
of Aff $S(G,u)$ and $\phi(G^+)$ is dense in Aff $S(G,u)^+$
(Corollary 13.7 and Theorem 14.14 of \cite{G}).

To obtain a consistent classification one imposes restrictions either
on the state space $S(G,u)$ or on the group $G$ itself. One reasonable
assumption is that $S(G,u)$ is finite-dimensional. In this case the
space Aff $(G,u)$ is isomorphic to ${\Bbb R}^k$ for some $k$ and
the classification of $(G,u)$ is reduced to classification of all dense 
subgroups of ${\Bbb R}^n$, cf Goodearl \cite{G}.

We take the second option making restrictions on $G$. For the topological
applications $G$ is the Betti group of a CW complex (Section 4).
Such groups are known to be free abelian of a finite rank $k\ge 2$.
Surprisingly, this assumption implies that the state space $S(G,u)$
is finite-dimensional. (Thus we come back to the first case.)
\begin{lem} \label{lm6}
{\bf (D.~Handelman, unpublished)}
Let $(G,u)$ be a simple interpolation group where $G$ is a free
abelian group of finite rank $k\ge 2$. Then the dimension of
$S(G,u)$ is at most $k-2$.
\end{lem}
{\it Proof.} Denote by $\{x_1,\dots,x_k\}$ a basis for $G$ as
a free abelian group. Let $\{p_1,\dots,p_k\}$ be the set of projections
on the corresponding coordinates (i.e. a dual basis). This set is
a basis for the real vector space $Hom_{\Bbb Z}(G,{\Bbb R})$.

Let to the contrary, the dimension of $S(G,u)$ is greater or equal to
$k-1$. Then $S(G,u)$ contains $n$ affinely independent states which
we denote by $s_1,\dots,s_k$. It is not hard to see that these states
are also linearly independent. Hence the vectors $\{s_1,\dots,s_k\}$
span $Hom_{\Bbb Z}(G,{\Bbb R})$. Therefore, there exists a nondegenerate
$k\times k$ matrix $(\alpha_{ij})$ with real entries such that 
$p_i=\sum_{j=1}^k\alpha_{ij}s_j$. 

Fix a number $m>k~max|\alpha_{ij}|$. Since $G$ is a free abelian
group of rank $k>1$, $G$ cannot be cyclic. Then there exists an
element $y\in G$ such that $y>0$ and $my<u$ ($u$ is the order unit).
Therefore $0\le s_j(y)\le 1/m$ for any $j$. Thus for all $i$
the following inequality is true:
\displaymath
|p_i(y)|\le\sum_{j=1}^k|\alpha_{ij}||s_j(y)|\le
\sum_{j=1}^k{|\alpha_{ij}|\over m}<\sum_{j=1}^k{1\over k}=1.
\enddisplaymath
Since $p_i(G)={\Bbb Z}$ for every $i$ we conclude that $p_i(y)=0$
and therefore $y=0$. This is a contradiction which proves
Lemma \ref{lm6}. 
$\square$
\begin{cor}\label{cr1}
Let $(G,u)$ be as in Lemma \ref{lm6}. Then $(G,u)$ is totally
ordered if and only if $S(G,u)$ consists of a unique state.
In this case the positive cone $G^+$ is defined by a hyperlane
passing through the origin of the space ${\Bbb R}^k$.
\end{cor}
{\it Proof.} If $(G,u)$ is totally ordered then $S(G,u)$
consists of a point (Corollary 4.17 of \cite{G}). If there
exists a unique state $s: (G,u)\to ({\Bbb R},1)$ then
$s$ is given by a linear function $f: {\Bbb R}^k\to {\Bbb R}$. 
The set $f^{-1}(0)$ is a hyperplane in ${\Bbb R}^k$ which splits
${\Bbb R}^k$ into two parts. Placing $G\cong {\Bbb Z}^k$ in
${\Bbb R}^k$ one gets a total order by the formula (\ref{eq1}).   
$\square$

\section{Main result}    
The space of totally ordered abelian groups of rank $k$ can be made to a
bundle space, $E_X$, over a CW complex $X$. Surprisingly enough,
the classification of $E_X$ reduces to  classification of certain
subbundles of a vector bundle $\xi$ of rank $k$ over $X$. An invariant 
responsible for such a classification turns out to be the Stiefel-Whitney 
class $w_1(\xi)$. (In \cite{Nik} bundles $E_X$ were called the characteristic 
fibrations.)

Let $(G,u)$ be a totally ordered simple abelian group of rank $k$.
By Corollary \ref{cr1} $(G,u)$ is defined by a hyperplane in the
Euclidean space ${\Bbb R}^k$. Let us denote by ${\cal S}_k(G,u)$ the
set of all such groups of a fixed rank $k$.
If $V_m({\Bbb R}^k)$ and $G_m({\Bbb R}^k)$ are the
Stiefel and Grassmann varieties, respectively, then by
$SG_m({\Bbb R}^k)$ we understand a subvariety of $G_m({\Bbb R}^k)$ 
consisting of $m$-dimensional {\it oriented} planes in ${\Bbb R}$. 
The manifold $SG_m({\Bbb R}^k)$ is known to be a double cover of 
the Grassmann manifold $G_m({\Bbb R}^k)$. There exists a natural function
\begin{equation}\label{eq14}
f: V_m({\Bbb R}^k)\longrightarrow SG_m({\Bbb R}^k),
\end{equation}
which identifies the $m$-frame $(v_1,\dots,v_m)$ with the linear
space $\langle v_1,\dots,v_m\rangle$ that it spans and which discerns the `left' and
`right' orientation of $m$-frames. Mapping $f$ is surjective since every
$m$-plane has a frame. The topology of $SG_m({\Bbb R}^k)$ coincides
with the topology induced by the mapping $f$. 

Let $m=k-1$ be a hyperplane (with orientation) in ${\Bbb R}^k$.   
By Corollary \ref{cr1} and discussion of Section 1, to an everywhere
dense
\footnote{In effect, this set is of `full measure' in the space
of ordered abelian groups of rank $k$ (not necessarily simple).
The residue set consists of groups with the order ideals. Roughly
speaking, such groups correspond to `rational directions' of the
hyperplane. Even if the direction is `irrational' it may happen
that the corresponding group is not totally ordered. However, such
a possibility is also exceptional (i.e. belongs to a residue set). 
This was shown in our work {\it Geometry of the Bratteli diagrams}, 
Preprint FI-NP2000-002, Fields Inst., 2000.}
subset, $\Omega$, of such hyperplanes there correspond a set of 
totally ordered simple abelian groups of rank $k$. Taking a closure
of $\Omega$ in $SG_{k-1}({\Bbb R}^k)$, one gets an {\it induced topology}
on the set ${\cal S}_k(G,u)$. The closure of ${\cal S}_k(G,u)$ is
denoted by $\bar{\cal S}_k(G,u)$. 
If $X$ is a CW complex of dimension $n$, then to every point of $X$
one can associate a dimension group $(G,u)\in \bar{\cal S}_k(G,u)$.
In other words, we have a continuous map
\begin{equation}\label{eq15}
X\longrightarrow \bar{\cal S}_k(G,u).
\end{equation}
Map (\ref{eq15}) is  cross-section of a locally trivial
bundle, which we denote by $E_X$. Globally there exists
a variety of topologically distinct bundles $E_X$, since
the pieces $U\times\bar{\cal S}_k(G,u)$ can be glued
together in a variety of ways. The following theorem establishes  
classification of such bundles.  
\begin{thm}\label{thm1}
Let $X$ be a finite CW complex. Then the topological class of
bundle $E_X$ is determined by the elements of group
$H^1(X; {\Bbb Z}_2)$. 
\end{thm}

\section{Proof of Theorem 1}    
The idea of the proof is a stepwise reduction of $E_X$ to
a subbundle of the vector bundle over $X$ whose topological
class is determined by the Stiefel-Whitney invariants. Of
course, there is no hope whatsoever that the obtained invariants
are sufficient, i.e. that they completely discern the topological classes
of $E_X$. But if $E_X$ and $E_X'$ has different set of invariants,
then they are topologically distinct.

Let $X, V_m({\Bbb R}^k)$ and $SG_m({\Bbb R}^k)$ be the manifolds
introduced in Section 2. Consider the following commutative diagram:

\bigskip
\begin{picture}(300,100)(70,0)
\put(150,5){$V_m({\Bbb R}^k)$}
\put(230,17){$f$}
\put(300,45){$\sigma$}
\put(210,60){$p$}
\put(190,8){\vector(3,0){73}}
\put(285,75){$X$}
\put(290,67){\vector(0,-1){45}}
\multiput(190,30)(16,8){4}{\line(2,1){10}}\put(254,62){\vector(2,1){13}}
\put(275,5){$SG_m({\Bbb R}^k)$}
\end{picture}

\bigskip\noindent
where $f$ is surjective and $\sigma$ is a cross-section of the
bundle $E_X$, see (\ref{eq15}). The map $p$ is defined so that
for any point $x\in SG_m({\Bbb R}^k)$ $p(f^{-1}(x))=y\in X$ such that 
$\sigma(y)=x$. In other words, $p$ sends the $m$-frames which define an 
oriented $m$-plane to the point of $X$ to which this plane is attached. 
Of course, map $p$ is not injective since to each $m$-plane in
$SG_m({\Bbb R}^k)$ one can relate a bunch of $m$-frames that span
this plane. To make $p$ injective, one must choose one $m$-frame
in every $m$-plane. The following lemma will be helpful.
\begin{lem}\label{lm11}
To every $m$-plane $\alpha\in SG_m({\Bbb R}^k)$ one can
relate a (left) $m$-frame of the space $V_m({\Bbb R}^k)$ 
defined uniquely by $\alpha$.
\end{lem}
{\it Proof.}
Denote by $O_n$ the orthogonal group of the real Euclidean
space ${\Bbb R}^n$. The Stiefel manifold $V_m({\Bbb R}^k)$
can be identified with the factor space $O_k/O_{k-m}$, see 
e.g. N.~Steenrod, {\it The topology
of fibre bundles,} Princeton, N.J., 1951, \S 7.5-7.10. In
this setting, the Grassmann manifold $G_m({\Bbb R}^k)
\cong O_k/(O_m\times O_{k-m})$ what means that the Stiefel
manifold is a bundle over $G_m({\Bbb R}^k)$ with the fibre
$O_m$.

Let us construct a standard $m$-frame in the fibre $O_m$. For 
that let us decompose $O_m$ into a sequence of inclusions:
\begin{equation}\label{eq16}
O_1\subset O_2\subset\dots\subset O_{m-1}\subset O_m.
\end{equation}
Fix a unit vector $v_1\in {\Bbb R}^1$. The orthogonal
group $O_1$ is cyclic of order 2, so it flips $v_1$
relatively the origin. We choose $v_2\in {\Bbb R}^2$
so that $O_2\supset O_1$ keeps invariant $v_1$ and complements
$v_1$ to a {\it left} orthonormal 2-frame. The procedure
can be prolonged by induction so that it stops with the
unit vector $v_m$. Thus, we have constructed a left $m$-frame
which is invariant under the action of group $O_{m+1}$. 
This frame we call {\it standard} and it was
obtained from sequence (\ref{eq16}) by a Gram-Schmidt process. 
In the standard frame matrix $O_m$ has a Jordan normal form which is 
block-diagonal. Lemma is proved.  
$\square$

\bigskip\noindent
Put $m=k-1$ in Lemma \ref{lm11} and apply the commutative
diagram that we discussed earlier. Equation (\ref{eq15}) 
in this case takes the form:
\begin{equation}\label{eq17}
X\longrightarrow V_{k-1}({\Bbb R}^k),
\end{equation}
where $V_{k-1}({\Bbb R}^k)$ is a Stiefel variety of $(k-1)$-frames
in ${\Bbb R}^k$. Map (\ref{eq17}) is the cross-section to a rank $k$ vector
bundle, ${\cal E}_X$, which admits in total $k-1$ linearly independent cross-sections 
over $X$. Two (topologically) distinct bundles ${\cal E}_X\ne {\cal E}_X'$ 
will result in two distinct bundles $E_X\ne E_X'$. The following lemma
is basic for the classification of bundles ${\cal E}_X$.   
\begin{lem}\label{lm12}
Let ${\cal E}_X$ and ${\cal E}_X'$ be two rank $k$ vector bundles
over a CW complex $X$ which admit $k-1$ linearly independent cross-sections.
If the first Stiefel-Whitney classes $w_1({\cal E}_X)$ and $w_1({\cal E}_X')$
are distinct, then the bundles ${\cal E}_X$ and ${\cal E}_X'$ are
topologically distinct.
\end{lem}
{\it Proof.} Since the vector bundle ${\cal E}_X$ has $k-1$ linearly independent
cross-section, the higher Stiefel-Whitney
classes of ${\cal E}_X$ vanish (Lemma \ref{lm9} of Appendix): 
\begin{equation}\label{eq18}
w_2({\cal E}_X)=w_3({\cal E}_X)=\dots= w_k({\cal E}_X)=0.
\end{equation}
Let us show that the remaining class $w_1({\cal E}_X)$ is
non-trivial. Indeed, denote by $\eta$ a rank $1$ vector
bundle over $X$ such that ${\cal E}_X\oplus\eta=\varepsilon$,
where $\varepsilon$ is a trival vector bundle of rank $k+1$.
By equations (\ref{eq18}) the total Stiefel-Whitney classes
of ${\cal E}_X$  and $\eta$ are as follows:
\begin{equation}\label{eq19}
w({\cal E}_X)=1+w_1({\cal E}_X),\qquad w(\eta)=1+w_1(\eta).
\end{equation}
Since the Whitney sum ${\cal E}_X\oplus\eta$ is trivial,
one can apply formula (\ref{eq9}) (cf. Appendix) so that 
$w(\eta)=\bar w({\cal E}_X)$,
where $\bar w({\cal E}_X)=1+w_1({\cal E}_X)+\dots$, see
Milnor-Stasheff \cite{MS}. Thus we have 
\begin{equation}\label{eq20}
w_1({\cal E}_X)=w_1(\eta)\ne 0,
\end{equation}
since there exists a plenty of non-trivial line bundles over
$X$. 

(Note that the same result can be proved by the methods of
obstruction theory, see \S 5.2. Namely, the Stiefel lemma
(Lemma \ref{lm10}) says that for $m=k-1$ the map
$f: X\to V_{k-1}({\Bbb R}^k)$ can always be continuously 
extended to a 1-dimensional skeleton, $K^1$, of $X$. This is 
exactly what equation (\ref{eq20}) suggests.)

Lemma is proved.
$\square$

\bigskip
To finish the proof of Theorem \ref{thm1} it remains to 
apply Lemma \ref{lm12} and axiom (ii) of the Stiefel-Whitney
classes, see \S 5.2. 

Theorem follows.
$\square$

\section{Characteristic classes of regular foliations}    
The bundles of ordered abelian groups  were of modest interest  
had  they no application outside noncommutative algebra. Fortunately, they
arise in the study of characteristic classes
of regular foliations. In particular, Theorem \ref{thm1} defines
the characteristic classes of such foliations. This section is reserved
for the discussion of some geometric  consequences of Theorem \ref{thm1}. For a
full treatment, see \cite{Nik}.

\subsection{Regular foliations}
Let us fix the following notation:

\bigskip
\begin{tabular}{cl}
$X$                 & $n$-dimensional hyperbolic manifold;\cr
$k$                 & first Betti number of $X$;\cr
${\Bbb H}^n$        & hyperbolic $n$-space $\{x\in {\Bbb R}^n|x_n>0\}$ endowed 
                      with metric\cr
                    & $ds=|dx|/x_n$;\cr
$\partial{\Bbb H}^n$& absolute $\{x\in {\Bbb H}^n~|~x_n=0\}$;\cr
$G$                 & M\"{o}bius group, i.e a discrete subgroup of 
                      $SL(2,{\Bbb C})/\pm I$;\cr
$X={\Bbb H}^n/G$    & representation of $X$ by actions of $G$ at the universal\cr
                    & cover ${\Bbb H}^n$;\cr
$\tau$              & tessellation of ${\Bbb H}^n$ by the fundamental domains\cr
                    & of group $G$;\cr
${\Bbb Z}^k\subseteq\tau$ & abelianized fundamental group of $X$, 
                            see details in \cite{Nik}.
\end{tabular}

\bigskip
Let $\cal F$ be a foliation on $X$ of the codimension 1. Such
foliations always exist on $3$-dimensional manifolds and in the higher 
dimensions good examples are known  (Lawson \cite{Law}).
Foliation $\cal F$ on $X$ will be called {\it regular} if 
it has no Reeb components and either of the following
conditions is satisfied: 

\medskip
(i) all leaves of $\cal F$  are  compact, i.e. the limit set
of every leaf is the leaf itself;

\smallskip
(ii) all leaves of $\cal F$  are dense in $X$,
i.e. the limit set of every leaf is the whole $X$, or a Cantor
set in $X$.

\bigskip\noindent
(One can usefully think of regular foliations as a higher-dimensional
analog of the Kronecker foliations on the two-dimensional torus.
These can either be ``rational'', i.e. with all leaves compact,
or ``irrational'', so that each leaf is everywhere dense on the
torus. An exceptional ``Denjoy'' case adds to the picture when
the leaves tend to a Cantor set.)

\bigskip\noindent
\underline{Asymptotic of regular foliations}.
Let $\cal F$ be a regular foliation. Denote by $\tilde {\cal F}$ 
the preimage of $\cal F$ on the universal cover
${\Bbb H}^n$. The asymptotic properties of leaves of 
$\tilde {\cal F}$ were studied by S.~Fenley, D.~Gabai and
W.~Thurston. They proved that any leaf $\tilde l$ of such foliation
is a quasi-isometric immersion into ${\Bbb H}^n$
and the limit set of $\tilde l$ is either a `circle' $S^{n-2}$
of the absolute or the whole absolute.
(See S.~Fenley, {\it Foliations with good geometry}, J. Amer. Math. Soc. 12, 1999,
619-676.)

\subsection{$({\Bbb Z}^k)^+$}
Every regular foliation generates a total order on the abelianized
fundamental group ${\Bbb Z}^k$ of $X$. For that take a leaf $\tilde l$ which 
has the limit `circle' $S^{n-2}$ at the absolute. This circle defines
a geodesic hemi-sphere $S^{n-1}$ in the hyperbolic $n$-space
which is either rational or irrational depending on property (i)
or (ii) of $\cal F$. A {\it positive cone} of ${\Bbb Z}^k\subseteq\tau$
is defined as
\begin{equation}\label{eq25}
({\Bbb Z}^k)^+=\{z\in {\Bbb Z}^k ~|~z\in~Int~S^{n-1}\},
\end{equation}
where $Int~S^{n-1}$ denotes the interior points of ${\Bbb H}^n$ separated
by the hemi-sphere $S^{n-1}$. The group $({\Bbb Z}^k,({\Bbb Z}^k)^+)$  is
called  an {\it ordered abelian group}  associated to $\cal F$.
This group is simple if $\cal F$  has property (ii)
and has non-trivial order ideals otherwise. A link between 
$({\Bbb Z}^k,({\Bbb Z}^k)^+)$ and the ring structure of a foliation $C^*$-algebra
introduced by A.~Connes is  given by the following lemma.   
\begin{lem}
The $K_0$-group of the $C^*$-algebra associated
to a regular foliation $\cal F$ is order-isomorphic to 
the totally ordered abelian group $({\Bbb Z}^k,({\Bbb Z}^k)^+)$.
\end{lem}
{\it Proof.} See \cite{Nik1}-\cite{Nik}.
$\square$

\subsection{$E_X$}
Given a regular foliation $\cal F$ on the manifold $X$
one can look at the {\it integrable distribution} which is
a set of planes tangent to $\cal F$ at every point of $X$.
We allow to rotate every plane of this distribution
with a single restriction that the newly obtained distributions
are integrable and correspond to regular foliations. The set of regular 
foliations obtained by a continuous rotation from $\cal F$ we call a
{\it homotopy class} of $\cal F$.  The elements of this homotopy class
can be put in bijection with the points of manifold $X$ as follows.

Consider a restriction of the regular foliation $\cal F$
to the unit cylinder $Cyl_h=D^{n-1}\times [0,h]$, $||D||=1$, $h>0$
of height $h$.  This restriction looks like a level set of the
function $h=~Const$. Let us rotate the distribution in $Cyl_h$ so that all the
tangent planes remain parallel. In this way, the rotation
is parametrized by a vector $v\in {\Bbb R}^n$ normal to
all planes. We set $||v||=h$. In this way  local deformations are 
parametrized by vectors $v$ whose endpoints $x\in {\Bbb R}^n$ we identify with
the points of manifold $X$ covered by a local chart in $0$. 
Gluing together the local homotopy deformations one gets a global one 
by identifying deformations on the overlapping domains.

(The reader should firmly understand that even local homotopy
deformations change dramatically the structure of foliations.
For example, an irrational foliation on the torus can be made rational,
and vice versa, by a local modification near any point of torus;
see also examples of Ch.~Pugh and others.) 
\begin{dfn}
By a bundle $E_X$ over manifold $X$ one understands a continuous
field of totally ordered abelian groups over $X$. Given $x\in X$
we relate a group $({\Bbb Z}^k,({\Bbb Z}^k)^+)$ if and only this
group corresponds to a regular foliation parametrized by $x$. 
In this way, the homotopy (characteristic) classes of regular
foliations correspond to topologically distinct types of $E_X$.  
\end{dfn}
\begin{cor}\label{cor2}
 {\bf (Characteristic classes of regular foliations)}
Let $X$ be a hyperbolic manifold whose first Betti number does not vanish. 
Then the characteristic classes of regular foliations on $X$
coincide with the elements of the first cohomology group $H^1(X;{\Bbb Z}_2)$
of $X$.
\end{cor}
{\it Proof.} The manifold $X$ admits a regular triangulation to a
CW complex for which we keep the same notation. Since $k\ne 0$,
one  applies Theorem \ref{thm1}. Corollary follows. 
$\square$

\section{Appendix: Stiefel-Whitney invariants}    
The characteristic classes were created by H.~Hopf, E.~Stiefel
and H.~Whitney to extend the Euler characteristic of a manifold.
Namely, the number and type of singularities of a vector field
is an invariant of the manifold (Euler number). Similarly, the
singularities of $n$-tuples of linearly independent vector fields
is an invariant of the vector bundle over the manifold (characteristic
class). For an introduction in the area we recommend the monograph
of Milnor and Stasheff \cite{MS} and the original papers of Stiefel 
\cite{Sti} and Whitney \cite{Whi}.

\subsection{Stiefel-Whitney class of a vector bundle}
We fix the following notation:

\bigskip
\begin{tabular}{cl}
$X$          & finite-dimensional CW complex;\cr
$\xi$        &vector bundle of rank $k$ over $X$ with total space $E$;\cr
$\xi\oplus\eta$     & Whitney sum of the vector bundles $\xi$ and $\eta$;\cr
$V_m({\Bbb R}^k)$   & Stiefel manifold of the orthogonal $m$-frames \cr
                    & in the Euclidean space ${\Bbb R}^k$;\cr
$G_m({\Bbb R}^k)$   & Grassmann manifold of the $m$-dimensional linear\cr
                    & subspaces in ${\Bbb R}^k$;\cr
$H^m(X;G)$          & singular cohomology group of $X$ of order $m$\cr
                    & with coefficients in the ring $G$. 
\end{tabular}

\bigskip\noindent
Let $G\cong {\Bbb Z}_2$ be the cyclic group of order $2$. By a
{\it Stiefel-Whitney class} of the vector bundle $\xi$ one understands
a sequence of the singular cohomology classes $\{w_i\in H^i(X;{\Bbb Z}_2)
~|~i=0,1,2,\dots\}$ such that the following axioms are satisfied:

\medskip
(i) $w_0(\xi)$ is a unit element of the group $H^0(X;{\Bbb Z}_2)$ and 
$w_i(\xi)=0$ for $i>k$ where $k$ is the rank of vector bundle $\xi$; 

\smallskip
(ii) if a bundle map $\xi\to\eta$ covers the map $X(\xi)\to X(\eta)$
then $w_i(\xi)=f^*w_i(\eta)$; 
\footnote{In particular, if $\xi$ and $\eta$ are topologically 
equivalent then $w_i(\xi)=w_i(\eta)$.}

\smallskip
(iii) cohomology class of the Whitney sum is calculated by the formula
\begin{equation}\label{eq6}
w_i(\xi\oplus\eta)=\sum_{j=1}^i w_j(\xi)\bigcup w_{i-j}(\eta);
\end{equation}

(iv) Stiefel-Whitney class $w_1$ of the line bundle over projective line is 
non-trivial.

\bigskip\noindent
The characteristic classes satisfying axioms (i)-(iv) do exist. A proof
of this fact can be found in Milnor-Stasheff \cite{MS}, \S 8. The rest 
of this subsection is devoted to some implications of axioms (i)-(iv).
\begin{lem}\label{lm7}
Let $\varepsilon$ be a trivial vector bundle. Then $w_i(\varepsilon)=0$
whenever $i>0$. 
\end{lem}
{\it Proof.} There exists a bundle map from $\varepsilon$ 
to a bundle $\eta$ whose base $X$ is a singleton. Lemma follows
from axiom (ii).  
$\square$
\begin{lem}\label{lm8}
Let $\varepsilon$ be a trivial vector bundle. Then $w_i(\varepsilon\oplus\eta)=w_i(\eta)$.
\end{lem}
{\it Proof.} This follows from Lemma \ref{lm7} and formula (\ref{eq6}).
$\square$
\begin{lem}\label{lm9}
Let $\xi$ be vector bundle of rank $k$. If $\xi$ admits $m$ cross-sections
$s: X\to E$ which define an $m$-frame at every point $p\in X$, then
\begin{equation}\label{eq7}
w_{k-m+1}(\xi)=w_{k-m+2}(\xi)=\dots = w_k(\xi)=0.
\end{equation}
In particular, if $\xi$ has a cross-section, then $w_k(\xi)=0$.
\end{lem}
{\it Proof.} Let $\xi$ have $m$ linearly independent cross-sections.
Then one can decompose $\xi=\varepsilon\oplus\varepsilon'$, where 
$\varepsilon$ is a trivial bundle of rank $m$ and $\varepsilon'$
is an orthogonal vector bundle of rank $k-m$. By Lemma \ref{lm8}
$w_i(\varepsilon\oplus\varepsilon')=w_i(\varepsilon')$. On the other
hand, rank $\varepsilon'=k$.  By axiom (i) lemma follows.    
$\square$

\bigskip\noindent
\underline{Ring of invariants}.
Formula (\ref{eq6}) can be simplified by considering a ring 
$H^*(X;{\Bbb Z}_2)$\linebreak$\buildrel\rm def\over=H^0(X;{\Bbb Z}_2)\times 
H^1(X;{\Bbb Z}_2)\times\dots$ of formal infinite series $w_0+w_1+\dots$
endowed with the usual (polynomial) multiplication.
The sum of two classes of order $m$ denotes the union of two $m$-cocycles
and the product of order $m$ and $m'$ classes means an $mm'$-cocycle which 
is the topological product of $m$- and $m'$-cocycle. By a
{\it total Stiefel-Whitney class} of a vector bundle $\xi$ one
understands an element $w(\xi)=1+w_1(\xi)+w_2(\xi)+\dots$
of the ring $H^*(X;{\Bbb Z}_2)$. In this notation formula (\ref{eq6})
reduces to:
\begin{equation}\label{eq8}
w(\xi\oplus\eta)=w(\xi)w(\eta).
\end{equation}
Equation (\ref{eq8}) can be resolved uniquely relatively $w(\eta)$.
For this one introduces an `inverse' element $\bar w(\xi)=1+\bar w_1+
\bar w_2+\dots$, where $\bar w_1=w_1$, $\bar w_2=w_1^2+w_2$, etc,
see \cite{MS}. In particular, if the bundle $\xi\oplus\eta$ is
trivial, then 
\begin{equation}\label{eq9}
w(\eta)=\bar w(\xi).
\end{equation}

\subsection{Obstruction theory}
The Stiefel-Whitney classes appear as a topological `obstruction' 
to continuous extension of a field of $m$-frames from the 
skeleton $K^{\nu-1}$ to the skeleton $K^{\nu}$ of a CW complex $X$.
Such an extension is possible if and only if the cohomology classes
of certain cocycles in $X$ vanish. The group of coefficients in
this cohomology is either $\Bbb Z$ or ${\Bbb Z}_2$. Milnor and
Stasheff showed that these cohomology classes define (and are completely
defined by) the Stiefel-Whitney classes of vector bundles over $X$
(Milnor-Stasheff \cite{MS}, pp. 139-143).

\bigskip\noindent
\underline{Stiefel manifolds}. Let us denote by $V_m({\Bbb R}^k)$
a set of the orthonormal $m$-frames with centre at the origin of the
Euclidean space ${\Bbb R}^k$. The space $V_m({\Bbb R}^k)$ is endowed
with a topology in which two $m$-frames are `close' if and only
if the endpoints of the respective basis vectors on the unit sphere
$S^{k-1}$ are close. This topology turns $V_m({\Bbb R}^k)$ to a
compact manifold of dimension $m(2k-m-1)/2$  (\cite{Sti},
p.311) which is called a {\it Stiefel manifold}. 
The homotopy and homology groups of $V_m({\Bbb R}^k)$ are
as follows (Stiefel \cite{Sti}, \S 1.4):
\begin{eqnarray}
\pi_{k-m-1}(V_m({\Bbb R}^k)) &=& 0,\label{eq10}\\
H_{k-m}(V_m({\Bbb R}^k)) &=& \cases{{\Bbb Z},\quad\hbox{if}\quad k-m\quad \hbox{is even}\cr
                                  {\Bbb Z}_2, \quad\hbox{if}\quad k-m\quad\hbox{is odd.}\label{eq11}} 
\end{eqnarray}

Let $X$ be a CW complex of dimension $n$. Consider the representation of
$X=\cup_{\nu=0}^nK^{\nu}$  as a sum of $\nu$-dimensional skeletons.      
The extension of an $m$-frame field $f: K^{\nu-1}\to V_m({\Bbb R}^k)$ 
to the skeleton $K^{\nu}$ depends essentially on extension of $f$ from 
the boundary $\partial E^{\nu}$ of elementary cell $E^{\nu}$ to the
whole $E^{\nu}$. Since $\partial E^{\nu}$ is homeomorphic to the 
unit $(\nu-1)$-sphere, the problem is to extend the map
\begin{equation}\label{eq12}
f:S^{\nu-1}\longrightarrow V_m({\Bbb R}^k),
\end{equation}
to the interior of the ball bounded by $S^{\nu-1}$.  The answer to this
question depends drastically on $\nu$.
\begin{lem}\label{lm10} {\bf (Stiefel)}
Let $f$ be a continuous map (\ref{eq12}). Let $\alpha$ be a characteristic
number defined from the cohomology equation $f(S^{\nu-1})=\alpha Z_{\nu-1}$,
where $Z_{\nu-1}$ is a $(\nu-1)$-dimensional cycle in $V_m({\Bbb R}^k)$. 
If $0\le\nu\le k-m$ then $f$ can always be continuously extended to the interior of
$S^{\nu-1}$. If $k-m+1\le\nu\le n$ then $f$ can be continuously extended
to the interior of $S^{\nu-1}$ if and only $\alpha\in G$ is zero. In the later
case

\medskip
(i) $G\cong {\Bbb Z}$ if and only if $\nu$ is odd;

\smallskip
(ii) $G\cong {\Bbb Z}_2$ if and only if $\nu$ is even. 
\end{lem}
{\it Proof.} The first part of lemma follows from the equation (\ref{eq10})
since every $(k-m-1)$-dimensional sphere in $V_m({\Bbb R}^k)$ is simply connected.

To prove the second part of lemma it is enough to notice that the image
$f(S^{\nu-1})$ of the sphere $S^{\nu-1}$ is a cycle in $V_m({\Bbb R}^k)$ which is
proportional to a `basic' cycle $Z_{\nu-1}$. By formula (\ref{eq11})
the coefficient, $\alpha$, in this proportion belongs to ${\Bbb Z}$ 
if $\nu-1$ is even, or to ${\Bbb Z}_2$ if $\nu-1$ is odd. 
$\square$  
\begin{cor}\label{cr2}
For every value of $\nu$ such that $k-m+1\le\nu\le n$ there exists
an obstruction class 
\begin{equation}\label{eq13}
{\cal O}_{\nu}\in H^{\nu}(X;G),
\end{equation}
where $G$ is as indicated in Lemma \ref{lm10}. The field of $m$-frames on
the skeleton $K^{\nu-1}$ can be extended to $K^{\nu}$ if and only if 
${\cal O}_{\nu}$ vanishes. Moreover, the obstruction
classes ${\cal O}_{\nu}$ are equal
\footnote{If one identifies the respective classes by a natural homomorphism
$h: {\Bbb Z}\to {\Bbb Z}_2$, cf \cite{MS}.}
 to the Stiefel-Whitney classes $w_{\nu}(\xi)$ of a vector bundle of rank $k$ over $X$. 
\end{cor}
{\it Proof.} The proof of this statement coincides with those of Theorem 12.1
and remarks on p. 143 of Milnor-Stasheff \cite{MS}.  
$\square$

\section{Conclusions}
The assumption on regular foliation to be of codimension 1 can be
eventually dropped. There is no difficulty to present any foliation
of codimension $p$ as a transversal `product'
\begin{equation}\label{eq21}
{\cal F}_1\bigcap\dots\bigcap {\cal F}_p,
\end{equation}
where each of ${\cal F}_i, ~1\le i\le p$ is a regular codimension 1
foliation. The theory of characteristic classes of such foliations
can be constructed by similar methods. The following conjecture    
is probably true.
\begin{con}
Let $X$ be a manifold specified in Corollary \ref{cor2}.
Then the characteristic classes of codimension $p$ foliations
on $X$ coincide with the elements of the product cohomology
group
\begin{equation}\label{eq22}
H^1(X;{\Bbb Z}_2)\times\dots\times H^p(X; {\Bbb Z}_2).
\end{equation}
\end{con}

    
\end{document}